# Jointly optimization of passenger-route assignment and transfer incentivization scheme for a customized modular bus system*

Jianbiao Wang, Tomio Miwa, Takayuki Morikawa

*Abstract*—As an emerging travel mode, the modular vehicle system (MVS) is receiving increasingly attention. In particular, the operators could connect multiply modular vehicles as an assembled bus in response to the temporarily demand varies. Therefore, in this study, the MVS is adopted in the context of customized bus design to satisfy passengers' reserved travel demand. In addition, to increase the potential of the customized modular bus system, the transfer among buses can be considered to group the passengers with same or close destinations. However, the passengers will not actively transfer as it is viewed as the disutility. Thus, the appropriate incentivization should be provided. To this end, we jointly optimize the passenger-route assignment and the transfer incentivization scheme, in which the transfer demand is considered elastically under different incentivization. Then, the linearization approaches are adopted to transform the original nonlinear model into mixed integer linear programing model, which can be solved by the state-of-the-art solvers. The experiment on the small network reveals that the performance of the bus system with incentivization scheme is better than that without incentivization scheme, and the extent of such superiority depends on the total demand level in the system.

## I. Introduction

The customized bus is gaining popularity among the world. However, the reserved travel demand differs among different time, and the common operated bus with fixed capacity is unable to be adaptive to the dynamic demand [1]. Consequently, when the supply and demand is mismatched, the operated bus is likely to be too vacant or too crowd. As an alternative, the emerging modular vehicle system (MVS) provide a new solution. In particular, the modular vehicles can be connected/disconnected as an assembled bus, of which the capacity is changeable in response to dynamic passenger demands. Thus, this kind of flexible MVS can be adopted to design the customized bus system. In addition, to increase the potential of customized modular bus system, the transfer can be incorporated. Through transfer activity, the passengers with same or close destinations can be grouped together, which makes the bus system more efficient. However, the transfer activity costs travelers extra time and efforts. Such disutility impedes the willingness of transferring among buses. In order to offset such negative influence, the incentive strategies can be adopted to encourage passengers to transfer when necessary. It is noted that the actual number of transfers required by the operator should not exceed the number of passengers who are willing to take part into transfer. Typically, a higher incentive will increase the willingness of transfer and vice versa. Thus, the relationship between the incentives and the willingness of transfer should also be considered. Commonly, the logit model is widely adopted to formulate the willingness of transfer in response to the monetary incentive, thus in this study such elastic nonlinear relationship is adopted [2].

In all, this paper proposes a model to design a customized modular bus system, the task is to jointly optimize passenger-route assignment and transfer incentivization scheme. In addition, the travel route of each bus as well as the bus capacity will be determined together. Moreover, several linearization techniques are used to transform the nonlinear optimization into the mixed integer linear programming (MILP), then the problem can be readily solved by the commercial state-of-the-art solvers like GUROBI. Finally, an experiment is conducted to test the performance of the proposed system design.

## II. Literature review

The transfer increases the connectivity and flexibility of bus transit system. Thus, it needs to be well designed for the daily operation. To this end, several studies have been conducted. For example, chu et al. [3] optimized the bus timetables to synchronize the transfer with a mixed-integer linear programming model. Moreover, to offset the negative influence of transfer, gong et al. [4] design a transfer-based customized modular bus system in which the maximal number of transfers for each user is set as once. In addition, the school bus routing and scheduling problem with transfers was investigated by Bögl et al, [5] and the impact of transfers on the service level was analyzed. However, most of studies regarded the users always follow the suggestion of operator if the transfer is needed. Actually, the final transfer number is determined from both the operator side and user side. Even if the operator decides how many transfers are needed, the rational user would not always follow the suggestion because it will waste the time and cause the inconvenience. Thus, a more practical way is to incentivize the user to participate into the transfer if necessary.

Considering the incentivization is also widely studied in transportation researches. For example, Wang and Liao [6] developed user-based incentives schemes to involve user into the relocation of car-sharing system. In terms of bus system,

*This work was supported by JSPS KAKENHI Grant Number 19H02260 and Nagoya University Interdisciplinary Frontier Fellowship Grant Number JPMJFS2120. The authors gratefully acknowledge their support. (Corresponding author: Jianbiao Wang.)

Jianbiao Wang is with the Department of Civil and Environmental Engineering, Nagoya University, Nagoya 464-8603, Japan (e-mail: wang.jianbiao98@gmail.com)

Tomio Miwa is with the Institute of Materials and Systems for Sustainability, Nagoya University, Nagoya 464-8603, Japan (e-mail: miwa@nagoya-u.jp)

Takayuki Morikawa is with the Institute of Innovation for Future Society, Nagoya University, Nagoya 464-8603, Japan (e-mail: morikawa@nagoya-u.jp)

Luo et al. [7] explored the relationships between the incentive subsidy provided for bus company and their operated service level. Similarly, Vigren and Pyddoke [8] analyzed the impact of individual incentive on the ridership of public bus system. Although several studies on incentivizing bus trips have been done, few of aforementioned studies considered the incentive for transferring, which should be highlighted especially when designing such a flexible customized modular bus system.

## III. PROBLEM AND FORMULATION

### A. Problem description

In the customized bus platform, the operator first released the candidate stations they will serve and the timetables the assembled bus fleets will depart. It is noted that the candidate stations and the departure time are predetermined. Then, the passengers send their requests including boarding stations, alighting stations, and preferred departure time. Based on the preferred departure time, the total demand for each time (e.g. 8:15 am) can be obtained. It is noted that since the demand among different timetables are irrelevant, the bus system is designed separately for each service time. Based on the optimal design, the passenger will be informed with mobility solution, it includes whether they can be serviced or not, which assembled bus they should aboard, how the bus travel through the stations, whether they need transfer or not, how much is provided to incentivize the transfer, and where they need to transfer. Here, the incentive provided by the operator is in the discrete levels, which includes no incentive, low incentive, middle incentive, and high incentive levels. Each level corresponds to a designated monetary incentive cost (e.g. 1$). The detail workflow is shown in Fig. 1.

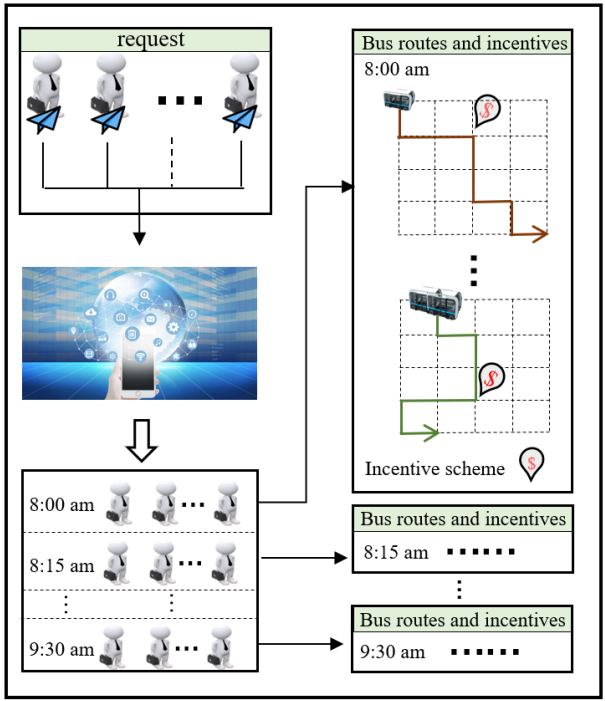

Fig. 1. The customized modular bus system design process

### B. Notions and formulations

Relevant notions are summarized in Table I, which is followed by the model formulation.

TABLE I. NOTION

| Notion | Definition |
|---|---|
| **Set** | |
| $\mathcal{F}$ | Set of stations, $i, j, m, n \in \mathcal{F}$ |
| $\mathcal{F}_i^+$ | Set of downstream stations of $i$ |
| $\mathcal{F}_i^-$ | Set of upstream stations of $i$ |
| $\mathcal{P}$ | Set of bus type $p \in \mathcal{P}$ |
| $\mathcal{K}$ | Set of bus route $k \in \mathcal{K}$ |
| $\mathcal{S}$ | Set of incentive strategy level $s \in \mathcal{S}$ |
| **Parameters** | |
| $c^e$ | Extra cost per unserved passenger. ($/c) |
| $c_p^o$ | Operation cost of type $p$ bus per kilometer ($/km) |
| $c^t$ | Value of travel time per passenger. ($/c) |
| $c^E$ | Total extra cost. ($) |
| $c^O$ | Total operation cost. ($) |
| $c^S$ | Total system cost. ($) |
| $c^T$ | Total passenger travel cost. ($) |
| $v_p$ | Capacity of type $p$ bus. (seats) |
| $d_{ij}$ | Distance between stations $i$ and $j$. ($km$) |
| $N$ | Total number of modular vehicle unit. |
| $q_{ij}$ | Passenger demand between OD pair $i$ and $j$. ($c$) |
| $t_{ij}$ | Travel time between stations $i$ and $j$. ($hour$) |
| $v$ | Capacity of a modular vehicle unit. (seats) |
| $\omega$ | weight coefficients in the objective function |
| **Decision variables** | |
| $x_{ijk}$ | 1 if link $(i, j)$ is served by bus $k$ and otherwise 0 |
| $y_{pk}$ | 1 if the type-p bus is used for route $k$ and otherwise 0 |
| $z_{ijmnk}$ | number of passengers between OD pair $(i,j)$ assigned to bus $k$ on link $(m, n)$ |
| $r_{ijmk}$ | number of passengers between OD $(i,j)$ transferring from bus $k$ to other vehicles at station $m$ |
| $p_{mks}$ | incentivize passengers to transfer from bus $k$ to other vehicles available at station $m$ with strategy $s$ |
| **Auxiliary variables for linearization** | |
| $b_{ijmk}$ | 1 if there are passengers between OD $(i,j)$ transferring from bus $k$ at station $m$, otherwise 0. |
| $w_{mkc}$ | 1 if user is willing to transfer from bus $k$ to other vehicles at station $m$ |
| $AU_{mkd}$ | the highest utility corresponds to the chosen alternative (i.e. transfer or not transfer from bus $k$ to other vehicles at station $m$) in scenario $d$, |
| $T_{ijmkd}$ | number of passengers between OD pair $(i,j)$ willing to transfer from bus k at station m in scenario $d$ |
| $R_{ijmks}$ | the incentive cost for user between OD pair $(i,j)$ transferring at station $m$ from bus $k$ arising from taking the $s^{th}$ incentivization strategy |

*1) Vehicle flow balance:* Each bus $k \in \mathcal{K}$ should travel through a complete route. Specifically, Constraint (1) and (2)

require that bus $k$ must start from and end at the virtual depot. Constrains (3) require that if bus $k$ enters station $j$, it must leave station $j$. Constraints (4) require that route $k$ only visits each station once by eliminating the subtours. For the convenience of the notation, define $s$ and $t$ as the start and end virtual depot.

$$\sum_{j\in\mathcal{F}} x_{sjk} = 1, \forall k \in \mathcal{K}. \quad (1)$$

$$\sum_{i\in\mathcal{F}} x_{itk} = 1, \forall k \in \mathcal{K}. \quad (2)$$

$$\sum_{i\in\mathcal{F}_j^-} x_{ijk} = \sum_{i\in\mathcal{F}_j^+} x_{jik}, \forall j \in \mathcal{F}\setminus\{s,t\}, k \in \mathcal{K}. \quad (3)$$

$$\sum_{i\in\not{f}} \sum_{j\in\not{f}_i^+} x_{ijk} \leq |\not{f}| - 1, \forall \not{f} \subset \mathcal{F}, 1 < |\not{f}| < |\mathcal{F}|, k \in \mathcal{K}. \quad (4)$$

*2) Passenger flow balance:* Constraints (5) guarantees that the number of serviced passengers between OD pair $(i,j)$ is less than or equal to the demand between $(i,j)$. Constraints (6) indicate the passenger flow balance at station $m$. Constraints (7) assures that, for OD pair $(i,j)$, the number of passengers arriving at destination $j$ must equal that departing from origin $i$. Constraints (8) requires that no passengers between OD pair $(i,j)$ will be assigned to bus $k$ on link $(m,n)$ unless bus $k$ travels on link $(m,n)$.

$$\sum_{k\in\mathcal{K}} \sum_{n\in\mathcal{F}_i^+} z_{ijink} \leq q_{ij}, \forall i \in \mathcal{F}, j \neq i \in \mathcal{F}. \quad (5)$$

$$\sum_{k\in\mathcal{K}} \sum_{n\in\mathcal{F}_m^-} z_{ijnmk} = \sum_{k\in\mathcal{K}} \sum_{n\in\mathcal{F}_m^+} z_{ijmnk},$$
$$\forall i \in \mathcal{F}, j \neq i \in \mathcal{F}, m \in \mathcal{F}\setminus\{i,j\} \quad (6)$$

$$\sum_{k\in\mathcal{K}} \sum_{m\in\mathcal{F}_j^-} z_{ijmjk} = \sum_{k\in\mathcal{K}} \sum_{n\in\mathcal{F}_i^+} z_{ijink}, \forall i \in \mathcal{F}, j \neq i \in \mathcal{F}. \quad (7)$$

$$z_{ijmnk} \leq q_{ij} x_{mnk},$$
$$\forall i \in \mathcal{F}, j \neq i \in \mathcal{F}, m \in \mathcal{F}\setminus\{j\}, n \in \mathcal{F}_m^+, k \in \mathcal{K}. \quad (8)$$

*3) Vehicle capacity constraint:* Constraints (9) assures that total passengers assigned to bus $k$ on each link is no more than the capacity of bus $k$. Constraints (10) indicate that bus $k$ is chosen as type $p$ assembled bus.

$$\sum_{i\in\mathcal{F}} \sum_{j\neq i\in\mathcal{F}} z_{ijmnk} \leq \sum_{p\in\mathcal{P}} y_{pk} pV, \forall m \in \mathcal{F}, n \in \mathcal{F}_m^+, k \in \mathcal{K}. \quad (9)$$

$$\sum_{p\in\mathcal{P}} y_{pk} = 1, \forall k \in \mathcal{K}. \quad (10)$$

*4) Passenger transfer:* Constraints (11) indicate that the number of passengers transferring from bus $k$ to other buses at station $m$ is the difference between the number of passengers arriving and leaving the station $m$ in bus $k$ if the former is larger than the latter, otherwise it equals to 0 meaning no transfers from bus $k$ at station $m$. Constraints (12) require that the number of transferred users from bus $k$ should be less or equal than the number of user willing to transfer arriving at station $m$. In (12), $\sum_{n\in\mathcal{F}_m^-} z_{ijnmk}$ is the passengers number arriving at station $m$, which is multiplied by the transfer willingness $P(ts_{mk})$ under the certain incentive strategy. The willingness $P(ts_{mk})$ is expressed as a logit form in (13). Specially, the base deterministic utility for non-transfer is set as $V_{mk1} = 0$, while the deterministic utility for transfer is set as $V_{mk2} = -cons + \sum_s ic_s p_{mks}$. The constant is set to a negative value to indicate the user's reluctance to participate in the transfer, but a stronger incentive strategy can increase the willingness of transfer. According to the discrete choice model, the total utility for transfer and non-transfer considering unobserved utility can be written as $U_{mkc} = V_{mkc} + \xi_{mkc}, c \in \{1,2\}$ respectively. The $\xi_{mkc}$ follow the Gumbel distribution. The constraints (14) indicate the operator must choose a strategy for bus $k$ at station $m$. The constraints (15) assures that there is no need to incentivize passengers if no one transfers from bus $k$ at station $m$, and in this case the no incentive strategy is chosen.

$$r_{ijmk} = max\{\sum_{n\in\mathcal{F}_m^-} z_{ijnmk} - \sum_{n\in\mathcal{F}_m^+} z_{ijmnk}, 0\},$$
$$\forall i \in \mathcal{F}, j \neq i \in \mathcal{F}, m \in \mathcal{F}\setminus\{i,j\}, k \in \mathcal{K}. \quad (11)$$

$$r_{ijmk} \leq P(ts_{mk}) \sum_{n\in\mathcal{F}_m^-} z_{ijnmk},$$
$$\forall i \in \mathcal{F}, j \neq i \in \mathcal{F}, m \in \mathcal{F}\setminus\{i,j\}, k \in \mathcal{K}. \quad (12)$$

$$P(ts_{mk}) = \frac{e^{-cons+\sum_s ic_s p_{mks}}}{e^{-cons+\sum_s ic_s p_{mks}} + e^0}, \forall m \in \mathcal{F}, k \in \mathcal{K}. \quad (13)$$

$$\sum_{s\in\mathcal{S}} p_{mks} = 1, m \in \mathcal{F}, \forall k \in \mathcal{K}. \quad (14)$$

$$\sum_{s\in\mathcal{S}\setminus\{1\}} p_{mks} \leq \sum_{i\in\mathcal{F}} \sum_{j\neq i\in\mathcal{F}} r_{ijmk}, \forall m \in \mathcal{F}, k \in \mathcal{K}. \quad (15)$$

*5) Maximum travel time and fleet size limit:* Constraints (16) suggest that the travel time of each route is no greater than a given threshold $\bar{T}$. Constraints (17) assures the number of modular vehicles for assembling is no more than the total modular vehicles the company has.

$$\sum_{i\in\mathcal{F}} \sum_{j\in\mathcal{F}_i^+} x_{ijk} t_{ij} \leq \bar{T}, \forall k \in \mathcal{K}. \quad (16)$$

$$\sum_{p\in\mathcal{P}} \sum_{k\in\mathcal{K}} p y_{pk} \leq N. \quad (17)$$

*6) Objective function:* The objective is to minimize total $C^S$, which is the sum of the operator cost $C^o$, the passenger travel cost $C^T$ and the unserved cost $C^E$,

$$\min_{x_{ijk}, y_{pk}, z_{ijmnk}, r_{ijmk}, p_{mks}} C^S = \omega_1 C^o + \omega_2 C^T + \omega_3 C^E, \quad (18)$$

where $\omega_1, \omega_2$ and $\omega_3$ are weights for each cost components. The operator cost $C^o$ includes the bus operation cost $C^{o1}$ and the passenger incentive cost $C^{o2}$. In particular, $C^{o1}$ is determined by the bus type and total length of each route, and $C^{o2}$ is determined by which incentive strategy is taken. Thus, $C^o$ can be formulated a

$$C^o = C^{o1} + C^{o2} = \sum_{k\in\mathcal{K}} \sum_{p\in\mathcal{P}} c_p^o y_{pk} \left(\sum_{i\in\mathcal{F}} \sum_{j\in\mathcal{F}_i^+} x_{ijk} l_{ij}\right) +$$
$$\sum_{i\in\mathcal{F}} \sum_{j\in\mathcal{F}\setminus\{i\}} \sum_{m\in\mathcal{F}\setminus\{i,j\}} \sum_{k\in\mathcal{K}} \sum_{s\in\mathcal{S}} ic_s p_{mks} \quad (19)$$

The passenger travel cost is defined as the product of the total travel time and the value of time $c^t$. In study [4], the transfer penalty is also considered, but here the penalty is ignored. In this study, the incentive is provided and the passengers choose to participate into transferring or not voluntarily. Thus, there is no need to consider the inconvenience cause by the transfer. The passenger travel cost can be formulated as

$$C^T = c^t \sum_{i \in \mathcal{F}} \sum_{j \in \mathcal{F} \setminus \{i\}} \left( \sum_{k \in \mathcal{K}} \sum_{m \in \mathcal{F} \setminus \{j\}} \sum_{n \in \mathcal{F}_m^+} z_{ijmnk} t_{mn} \right) \quad (20)$$

Finally, we define the extra cost for unserved passengers completing their travels by private cars. It is measured as the product of the number of unserved passengers and the unit extra cost $c^e$. Thus, the external cost is formulated as

$$C^E = c^e \sum_{i \in \mathcal{F}} \sum_{j \in \mathcal{F} \setminus \{i\}} \left( q_{ij} - \sum_{k \in \mathcal{K}} \sum_{n \in \mathcal{F}_i^+} z_{ijink} \right) \quad (21)$$

## IV. SOLUTION APPROACH

The formulated model (1)-(21) is a mixed integer non-linear programming model. To address the nonlinearity, we reformulate it into a mixed integer linear programming model (MILP). This new formulation allows us to solve the exact optimal solution with state-of-the-art commercial solvers. The non-linearity of model (1)-(21) comes from the fact that the number of transfers in (11) is a maximal function, the number of passenger willing to transfer in (12) is elastic, and the operational cost $C^{o1}$ (19) is a nonlinear function. Those three nonlinearities are linearized as follows.

*1) Linearization of (11):* The maximal expression in (11) is a non-linear expression. To linearize it, a binary auxiliary variable $b_{ijmk}$ is introduced. The linear form is given as:

$$r_{ijmk} \geq 0, \forall i \in \mathcal{F}, j \neq i \in \mathcal{F}, m \in \mathcal{F} \setminus \{i,j\}, k \in \mathcal{K}. \quad (22)$$

$$r_{ijmk} \geq \sum_{n \in \mathcal{F}_m^-} z_{ijnmk} - \sum_{n \in \mathcal{F}_m^+} z_{ijmnk},$$
$$\forall i \in \mathcal{F}, j \neq i \in \mathcal{F}, m \in \mathcal{F} \setminus \{i,j\}, k \in \mathcal{K}. \quad (23)$$

$$r_{ijmk} \leq \sum_{n \in \mathcal{F}_m^-} z_{ijnmk} - \sum_{n \in \mathcal{F}_m^+} z_{ijmnk} + M_{11}(1 - b_{ijmk}),$$
$$\forall i \in \mathcal{F}, j \neq i \in \mathcal{F}, m \in \mathcal{F} \setminus \{i,j\}, k \in \mathcal{K}. \quad (24)$$

$$r_{ijmk} \leq M_{11} b_{ijmk},$$
$$\forall i \in \mathcal{F}, j \neq i \in \mathcal{F}, m \in \mathcal{F} \setminus \{i,j\}, k \in \mathcal{K}. \quad (25)$$

where $M_{11}$ is a constant big value. If the value of $\sum_{n \in \mathcal{F}_m^-} z_{ijnmk} - \sum_{n \in \mathcal{F}_m^+} z_{ijmnk}$ is positive, the (22)-(25) will force the $b_{ijmk}$ to be 1 ($b_{ijmk} = 0$ is not a feasible solution), and consequently, the $r_{ijmk}$ reduces to $\sum_{n \in \mathcal{F}_m^-} z_{ijnmk} - \sum_{n \in \mathcal{F}_m^+} z_{ijmnk}$. Otherwise, the $r_{ijmk}$ is 0 with the indicator $b_{ijmk}$ equals to 0.

*2) Linearization of (12):* The next step is to linearize the elastic transfer demand in (12)-(13). To begin with, we introduce binary variables $w_{mkc}$, which equals to 1 if the alternative transfer is chosen and 0 otherwise. As the user under each case would exactly decide transfer or not transfer, we impose

$$\sum_{c \in \{1,2\}} w_{mkc} = 1, \forall m \in \mathcal{F}, k \in \mathcal{K}. \quad (26)$$

According to the utility-maximizing theory, the users are willing to transfer if they feel the utility of transfer is higher than that of non-transfer, i.e., $U_{mk1} > U_{mk2}$. However, the unobserved utilities $\xi_{mkc}$ in $U_{mkc}$ are uncertain value. To make utility function comparable, the simulation technique is adopted. The core idea of simulation technique is to draw the value of unobserved utility $\xi_{mkdc}$ for each $d \in D$ scenario. Therefore, the constraints (26) for are extended for each draw $d$:

$$\sum_{c \in \{1,2\}} w_{mkdc} = 1, \forall m \in \mathcal{F}, k \in \mathcal{K}, d \in D. \quad (27)$$

Once the uncertain terms are drawn, the utility function $U_{mkdc}$ is determined, and consequently the utility of alternatives can be compared. Thus, the value of $w_{mkdc}$ can be determined with $w_{mkdc} = 1$ for highest $U_{mkdc}$, and 0 for another. The continuous auxiliary variable $AU_{mkd} \in (-\infty, -\infty)$ is introduced to represent the highest utility corresponds to the chosen alternative (transfer or not transfer):

$$AU_{mkd} = \max_c \{U_{mkdc}\}, \forall m \in \mathcal{F}, k \in \mathcal{K}, d \in D. \quad (28)$$

Then, the linear formulation of (28) is given by:

$$U_{mkdc} \leq AU_{mkd}, \forall m \in \mathcal{F}, k \in \mathcal{K}, d \in D, c \in \{1,2\}. \quad (29)$$
$$AU_{mkd} \leq U_{mkdc} + M_{mkd}(1 - w_{mkdc}),$$
$$\forall m \in \mathcal{F}, k \in \mathcal{K}, d \in D, c \in \{1,2\}. \quad (30)$$

where $M_{mkd}$ is a big constant value. In (29)-(30), If the highest utility comes from transfer (i.e., $AU_{mkd} = U_{mkd1}$), then $w_{mkd1}$ will be forced to be 1 (i.e., user choose to transfer) and otherwise 0. As a consequence of the law of large numbers, the probability of user willing to transfer can be calculated by the sample average:

$$P(t|p_{mk}) = \frac{1}{D} \sum_{d=1}^{D} w_{mkd1}, \forall m \in \mathcal{F}, k \in \mathcal{K}. \quad (31)$$

By substituting (31) in (12), the upper bound for the number of transfers is

$$r_{ijmk} \leq \frac{1}{D} \sum_{d=1}^{D} w_{mkd1} \sum_{n \in \mathcal{F}_m^-} z_{ijnmk},$$
$$\forall i \in \mathcal{F}, j \neq i \in \mathcal{F}, m \in \mathcal{F} \setminus \{i,j\}, k \in \mathcal{K}. \quad (32)$$

However, the (32) is still non-linear since it involves the production of binary variable $w_{mkd1}$ and integer variable $z_{ijnmk}$. We then introduce the auxiliary variable $T_{ijmkd} \in \mathbb{R}^+ \cup \{0\}$, which represents the number of user willing to transfer if under scenario $d$. Moreover, we define a big constant value $M$. The (32) can be linearized as (33) together with additional constraints (34)-(36).

$$r_{ijmk} \leq \frac{1}{D} \sum_{d=1}^{D} T_{ijmkd},$$
$$\forall i \in \mathcal{F}, j \neq i \in \mathcal{F}, m \in \mathcal{F} \setminus \{i,j\}, k \in \mathcal{K}, \quad (33)$$

$$T_{ijmkd} \geq \sum_{n \in \mathcal{F}_m^-} z_{ijnmk} - M(1 - w_{mkd1}),$$
$$\forall i \in \mathcal{F}, j \neq i \in \mathcal{F}, m \in \mathcal{F} \setminus \{i,j\}, k \in \mathcal{K}, d \in D, \quad (34)$$
$$T_{ijmkd} \leq \sum_{n \in \mathcal{F}_m^-} z_{ijnmk} + M(1 - w_{mkd1}),$$
$$\forall i \in \mathcal{F}, j \neq i \in \mathcal{F}, m \in \mathcal{F} \setminus \{i,j\}, k \in \mathcal{K}, d \in D, \quad (35)$$
$$T_{ijmkd} \leq M w_{mkd1},$$
$$\forall i \in \mathcal{F}, j \neq i \in \mathcal{F}, m \in \mathcal{F} \setminus \{i,j\}, k \in \mathcal{K}, d \in D. \quad (36)$$

Specially, Eqs. (34)-(36) reduce to $T_{ijmkd} = \sum_{n \in \mathcal{F}_m^-} z_{ijnmk}$ when $w_{mkd1} = 1$ and to $T_{ijmkd} = 0$ when $w_{mkd1} = 0$. By summing up all the passengers willing to transfer in $D$ scenarios and divided by $D$, the expected number of user willingness to transfer under the given incentive strategy is calculated.

*3) Linearization of (19):* For (19), both the cost for operating bus $C^{o1}$ and incentivizing user $C^{o2}$ are non-linear function. Both the nonlinearity result from the production of binary variable and integer variable. First, the cost for operating bus is linearized as follows. We define an auxiliary variable $C_{kp}^o \in \mathbb{R}^+ \cup \{0\}$, which represents the operational cost of route $k$ arising from operating with type-$p$ bus. Define $M_p$ a big constant value, and then $C^{o1}$ can be reformulated as,

$$C^{o1} = \sum_{k \in \mathcal{K}'} \sum_{p \in P} C_{kp}^o, \quad (37)$$

$$C_{kp}^o \geq c_p^o \left( \sum_{i \in \mathcal{F}} \sum_{j \in \mathcal{F}_i^+} x_{ij} l_{ij} \right) - M_p(1 - y_{pk}),$$
$$\forall k \in \mathcal{K}, p \in \mathcal{P}. \quad (38)$$

$$C_{kp}^o \leq c_p^o \left( \sum_{i \in \mathcal{F}} \sum_{j \in \mathcal{F}_i^+} x_{ij} l_{ij} \right) + M_p(1 - y_{pk}), \forall k \in \mathcal{R}, p \in \mathcal{P}. \quad (39)$$

$$C_{kp}^o \leq M_p y_{pk}, \forall k \in \mathcal{K}, p \in \mathcal{P}. \quad (40)$$

It is observed that Eqs. (38)-(40) reduce to $C_{kp}^\circ = c_p^o \left( \sum_{i \in \mathcal{F}} \sum_{j \in \mathcal{F}_i^+} x_{ijk} k_{ij} \right)$ when $y_{pk} = 1$ and to $C_{kp}^\circ = 0$ when $y_{pk} = 0$, which is the same feasible region that $C^{o1}$ in (19) define. Similarly, the cost for incentivizing users $C^{o2}$ is linearized by introducing the auxiliary variable $R_{ijmks} \geq 0$. $R_{ijmks}$ represents the incentive cost for user between OD pair $(i,j)$ transferring at station $m$ from vehicle $k$ arising from taking $s^{th}$ incentivization strategy. A big constant $M$ is defined. Then,

$$C^{o2} = \sum_{i \in \mathcal{F}} \sum_{j \in \mathcal{F} \setminus \{i\}} \sum_{m \in \mathcal{F} \setminus \{i,j\}} \sum_{k \in \mathcal{K}} \sum_{s \in \mathcal{S}} R_{ijmks} \quad (41)$$

$$R_{ijmks} \geq ic_s r_{ijmk} - M(1 - p_{mks}),$$
$$\forall i \in \mathcal{F}, j \neq i \in \mathcal{F}, m \in \mathcal{F} \setminus \{i,j\}, k \in \mathcal{K}, s \in \mathcal{S} \quad (42)$$
$$R_{ijmks} \leq ic_s r_{ijmk} + M(1 - p_{mks}),$$
$$\forall i \in \mathcal{F}, j \neq i \in \mathcal{F}, m \in \mathcal{F} \setminus \{i,j\}, k \in \mathcal{K}, s \in \mathcal{S} \quad (43)$$
$$R_{ijmks} \leq M p_{mks},$$
$$\forall i \in \mathcal{F}, j \neq i \in \mathcal{F}, m \in \mathcal{F} \setminus \{i,j\}, k \in \mathcal{K}, s \in \mathcal{S} \quad (44)$$

Eqs. (42)-(44) reduce to $R_{ijmks} = ic_s r_{ijmk}$ when $p_{mks} = 1$ and to $R_{ijmks} = 0$ when $p_{mks} = 0$, which is the same feasible region that $C^{o2}$ in (19) define. With the above linearization, an equivalent linear programming (ELP) formulation of the original model can be obtained as follows.

**ELP:**
$$\min_{x_{ijk}, y_{pk}, z_{ijmnk}, r_{jmk}, p_{mks}} \omega_1 C^O + \omega_2 C^T + \omega_3 C^E. \quad (45)$$

**subject to:** Constraints (1)-(10), (14)-(17), (22)-(27), (29)-(30), (33)-(44).

## V. NUMERICAL EXPERIMENT

### A. Experiment setup

The following Fig. 2 shows case study network, which includes 7 bus stations and 9 bidirectional roads. In this study, three assembled bus routes are used for the service, i.e. $\mathcal{K} = 3$. The speed of each bus is assumed as 30 km/h. Moreover, there are three types of bus to choose, of which the capacity is 10, 20, and 30 respectively, i.e. $\mathcal{P} = \{1,2,3\}$, the associated operational costs ($/km) for the three types of bus are 0.35, 0.5 and 0.7, respectively [4]. In addition, to test the robustness of the proposed model, the level of demand between each OD pair is generated randomly from 1 to 5, and then the demand between each OD pair is calculated according to the total demand set for each scenario. As for the utility function for transferring, the constant value is set as -1. There are four incentive strategies for operators to choose, i.e., no incentive, low incentive, middle incentive, and high incentive, of which the incentive cost ($) is set as 0, 1, 2, and 3, respectively. The number of draws for linearizing (12) is set as 30, i.e., $D = 30$. By following the study of [4], the weight coefficient for the operational cost, passenger travel cost, and extra service cost in the objective function is set as 0.4, 0.4, and 0.2. As for the other parameters, the maximum route duration $\bar{T}$ in (16), fleet size in (17), value of time for passenger $c^t$ in (20), and unit extra cost $c^e$ in (21) is set as 26 min, 8, 6 $/passenger, and 3.7 $/hour, respectively [4].

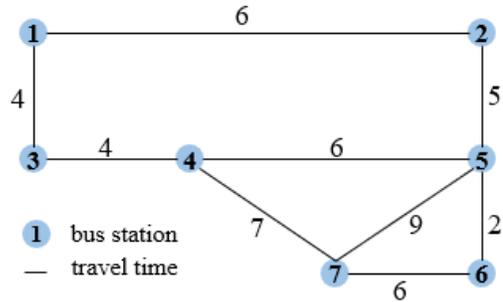

Fig. 2. bus stations and network

### B. Result discussion

We first illustrate a case design in detail in Table II with the total demand set as 100. Then the performance of the proposed design is compared with the system where the incentive for transfer is not provided (accordingly, no passenger is willing to transfer). By following the literature [4], five indicators in Table III are adopted for the comparison, which include the total travel distance of buses (TTD), the average in-vehicle travel time for each passenger (AIVTT), the transfer ratio defined as the percentage of transfers (TR), the service rate (SR) defined as the ratio of the number of passengers served to the total demand in the system, and the total system cost (TSC, i. e., the objective value).

TABLE II BUS OPERATIONAL NETWORK DESIGN

| Bus No. | Capacity | Path route | Transfer node | Transfer incentive | Transferred demand | Served demand |
|---------|----------|------------|---------------|--------------------|--------------------|---------------|
| 1 | 20 | 3→4→7→6→5→2 | 6 | high | 2 | 27 |
| 2 | 20 | 6→5→4→3→1→2 | n/a | n/a | n/a | 35 |
| 3 | 30 | 2→1→3→4→5 | 4 | high | 6 | 30 |

From Table II, the bus No.1 and No.2 are both assembled with two modular vehicles and both pass the 6 stations, while the bus No.3 assembled with three modular vehicles and passes 5 station. The passengers in bus No.1 and No.3 are highly incentivized to transfer at station 6 and 4, respectively. As a consequence, there are 2 and 6 passengers transferring at station 6 and 4, respectively. In all, the bus No.1, bus No.2, and bus No.3 served 27 passengers, 35 passengers, and 30 passengers, respectively.

The following Table III shows the comparison between the system with and without the incentivization scheme. As no one is willingness to transfer when the no incentivization is provided, the indicator TR of No incentive scheme is always 0. Thus, we compare the left 4 indicators. The result reveals that the 3 out of 4 indicators (i.e. TTD, SR, TSC) of Incentive scheme are better except for AIVTT. The increase of AIVTT under Incentive scheme results from the transfer activity. However, such a small increase of 0.04min can be ignored. Thus, it is concluded that the proposed incentive scheme could achieve better performance

TABLE III PERFORMANCE EVALUATION

|  | TTD | AIVTT | TR | SR | TSC |
|---|-----|-------|----|----|-----|
| Incentive | 33 km | 10.55 min | 9% | 92% | 47.8 $ |
| No incentive | 34 km | 10.51 min | 0% | 86% | 58.7 $ |

Then, the sensitivity analysis is conducted to test the performance between the *incentive* scheme and *no incentive* scheme under different demand level. To this end, we increase the total demand from 50 to 150. As the primary purpose of the bus system design is to satisfy travelers demand as much as possible, we mainly focus on the service rate in the sensitivity analysis. Besides, the transfer number under each demand level is also presented. The result is shown in Fig.3. For *incentive* scheme, the demand is fully satisfied until the total demand reaches 80. On the other hand, the service rate of *no incentive* scheme shows a continuous decrease trend. The service rate of *incentive* scheme is always higher than that of *no incentive* scheme, but the different between them is smaller with the increase of demand level. In addition, for the *incentive* scheme, the number of transfers remains at a stable level around 10 until the total demand exceed 100, then the transfer number begin to decrease.

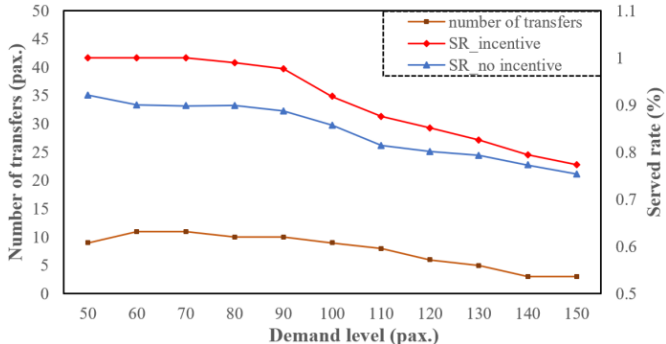

Fig.3. Sensitivity analysis

The trend of transfer number result from the fact that the room spared for the transfer decreases with the total demand increase. Moreover, the decrease of numbers of transfer indicates that the advantage of *incentive* scheme is gradually diminishing, which also explained why the difference between the service rate of two schemes becomes smaller. It reveals that to take full advantage of the *incentive* scheme, the total demand should be controlled within the certain level, otherwise the superiority of the *incentive* scheme is not significant. In addition, it is possible to maintain such superiority by increasing the number of modular vehicles or assembled bus fleet sizes, but it will increase the operation cost, such trade-off analysis is leaved for future analysis.

## VI. CONCLUSION

In this study, the customized modular bus system is designed by optimizing the passenger-route assignment and the transfer incentivization scheme jointly. Particularly, the transfer demand is considered elastically under different incentivization schemes. The nonlinear model is then transformed to the MILP, which can be readily solved by GORUBI. The experiment reveals that the performance of the bus system with incentivization scheme is better than that without incentivization scheme. In future, the efficient heuristic algorithms can be developed, with which the large case study in real world can be conducted.